\documentclass[a4paper,11pt]{amsart}

\usepackage{upref,amssymb}
\usepackage[leqno]{amsmath}
\usepackage[arrow,matrix]{xy}
\usepackage{pxfonts}
\usepackage{color}
\usepackage[margin=2cm]{geometry}

\newtheorem{prp}{Proposition}[section]

\newtheorem{thm}{Theorem}[section]

\theoremstyle{definition}
\newtheorem{rem}{Remark}

\newtheorem{lem}{Lemma}[section]
\theoremstyle{remark}

\title{Mod $2$ cohomology of $2$-local finite groups of low rank}
\author{Shizuo Kaji}
\thanks{\tiny The author was partially supported by the Grant-in-Aid for JSPS Fellows 182641.}


\address{Department of Mathematics \endgraf
                       Kyoto University \endgraf 
                       Kyoto 606-8502  \\  Japan}
\email{kaji@math.kyoto-u.ac.jp}

  \def\z{\mathbb{Z}}
\def\t{\mathcal{T}} 
\def\L{\mathcal{L}} \def\q{\mathbb{Q}}
\def\sig{\tilde{\sigma}_\phi}
\def\a2{\mathcal{A}_2}

\begin{document}
\maketitle

%%%%%%%%%%%%%%%%%%%%%%%%%%%%%%%%%%%%%%%%%%%%%%%%%
\begin{abstract}
  We determine the mod $2$ cohomology over the Steenrod algebra 
 $\a2$ of
 the classifying spaces of the free loop groups $LG$ for compact groups $G=Spin(7)$, $Spin(8)$, $Spin(9)$,
 and $F_4$. Then, we show that they are isomorphic as algebras over $\a2$
 to the mod $2$ cohomology of 
 the corresponding Chevalley groups of type $G(q)$, where $q$ is an odd prime power.
 In a similar manner, we compute the cohomology of the free loop space over $BDI(4)$ 
 and show that it is isomorphic to that of $BSol(q)$ as algebras over $\a2$.
 This note is a revised version of \cite{kaji}.
\end{abstract}

\section{Introduction}
For a based self-map $\phi: X\to X$ of a based topological space $X$, the \emph{homotopy fixed point space} (or the \emph{twisted loop space}) $\L_\phi X$ of $\phi$ is defined by the pullback
\begin{equation}\label{eq:sq}
\xymatrix{
\L_\phi X \ar[r] \ar[d] & X^{[0,1]} \ar[d]^{(ev_0,ev_1)} \\
X \ar[r]^(0.35){(Id,\phi)} & X\times X,
}
\end{equation}
where $ev_i: X^{[0,1]} \to X \ (i=0,1)$ are the maps assigning $\gamma(i)\in X$ to a path $\gamma\in  X^{[0,1]}$.
%In other words, $\L X$ is a space of all continuous maps $l$ from the interval $[0,1]$ to $X$ which satisfy $l(0)=\phi(l(1))$.
The free loop space $LX=\{\gamma: S^1\to X\}$ over $X$ is identified with the homotopy fixed point space of
the identity map of $X$.
Kuribayashi~\cite{Kuribayashi} introduced a map 
from the cohomology of $X$ to that of $LX$ with degree $-1$,
called the {\em module derivation},
to compute the cohomology
of $LX$ by the Eilenberg-Moore spectral sequence of the above pullback diagram.
Another interesting example of the 
homotopy fixed point space is that of 
an unstable Adams operation of
 the classifying space $BG$ of
 a compact Lie group $G$. 
Denote the Bousfield-Kan $2$-completion~\cite{Completion} 
of $X$ by $X^{\wedge}_2$.
Let $q$ be an odd prime power throughout this note.
There exists a self-map $\psi^q$ of $(BG)^{\wedge}_2$, called the unstable Adams operation
 of degree $q$ (\cite{Adams-operation}), 
 %which is an extension of the induced map of the $q$-th power map on the maximal torus.
 such that $\psi^q$ induces multiplication by $q^r$ on $H^{2r}((BG)^{\wedge}_2;\q)$.
The space $\L_{\psi^q} (BG)^{\wedge}_2$ is 
known to have the same mod $2$ homotopy type as the classifying space of
the finite Chevalley group of type $G(q)$~\cite{Etale}.
In particular, $H^*(G(q);\z/2)\cong H^*(\L_{\psi^q} (BG)^{\wedge}_2;\z/2)$.
Moreover, for the finite loop space $DI(4)$ at prime $2$ constructed by Dwyer and Wilkerson~\cite{Dwyer-Wilkerson},
Notbohm~\cite{BDI(4)} showed that there exists a self-map
$\psi^q$ of $BDI(4)$, also called the unstable Adams operation of
degree $q$, which induces multiplication by $q^r$ on $H^{2r}(BDI(4);\z_2^\wedge)\otimes \q$.
%which is an extension of the induced map of the $q$-th power map on the maximal torus.
Benson~\cite{Benson} defined the
classifying space $BSol(q)$ of an exotic $2$-local finite group 
 as $L_{\psi^q} BDI(4)$,   
which can be regarded as the ``classifying space'' of Solomon's non-existent finite
group~\cite{Solomon}.
%In this case, $\L_{\psi^q} BDI(4)$ is the classifying space $BSol(q)$ of the exotic $2$-local finite group defined by Benson (\cite{Benson}).

Kishimoto and Kono~\cite{Kishimoto-Kono}
generalized Kuribayashi's work to
calculate the cohomology of $\L_\phi X$.
Their result provides an efficient computational tool when the cohomology of $BG$ is a polynomial algebra.
As an example, they computed 
$H^*(LBG_2;\z/2)$ for the compact simple exceptional Lie group $G_2$ 
and showed that it is isomorphic to $H^*(BG_2(q);\z/2)$
as algebras over the Steenrod algebra.
 In this note, we carry out the computation of the mod $2$ cohomology of 
 $LBSpin(7), LBSpin(8), LBSpin(9), LBF_4$, and $LBDI(4)$
 as algebras over the Steenrod algebra $\a2$.
Furthermore, we show that they are isomorphic as algebras over the Steenrod algebra $\a2$
to those of the corresponding $\L_{\psi^q} (BG)^{\wedge}_2$.
The main theorem is as follows:
\begin{thm}\label{thm:1}
Let $q$ be an odd prime power.
   The following are isomorphisms of algebras over the Steenrod algebra:
   \begin{align*}
   H^*(LBSpin(n);\z/2) & \cong H^*(BSpin_n(q);\z/2)  \quad (n=7,8,9) \\
   H^*(LBF_4;\z/2) & \cong H^*(BF_4(q);\z/2)  \\
   H^*(LBDI(4);\z/2) & \cong H^*(BSol(q);\z/2).
  \end{align*}
\end{thm}
When $G$ is a compact Lie group, 
by the homotopy equivalence $BLG\simeq LBG$ (see for example \cite[\S 2]{Kuribayashi}),
the above theorem establishes cohomology isomorphisms between the classifying spaces of the loop groups $LG$
and those of the corresponding finite groups $G(q)$ for $G=Spin(7), Spin(8), Spin(9)$, and $F_4$.
The concrete presentations of the above algebras are given in Propositions 
\ref{prp:bspin7}, \ref{prp:bspin8}, \ref{prop:bspin9}, \ref{prp:bf4}, and \ref{prp:bdi4}.

{\it  Acknowledgment.}  We  would like to thank Professor Akira Kono 
  for various suggestions.     
  Also, we are very grateful to Anssi Lahtinen, who has corrected many errors in the previous version of this note.
  We have revised this note according to a detailed report kindly provided by him.

% {\it  Apology.} We made serious mistakes in the printed version of this note.
%   We would like to apologize sincerely here.

\section{Main tool}
Here, we summarize the results of Kishimoto and Kono~\cite{Kishimoto-Kono} that are necessary for our
purpose. Unless a coefficient ring is specified, $H^*(X)$ always means the mod $2$ cohomology of $X$.
%Let $\phi$ be a based self-map of a based space $X$.  The
%twisted loop space $\L X$ of $X$ is defined in the following
%pull-back diagram:
%\[ \xymatrix{
%  \L X \ar[r] \ar[d] & X^{[0,1]} \ar[d]^{e_0 \times e_1} \\
%  X \ar[r]^{1\times \phi} & X \times X }
%\]
%where $e_i \ (i=0,1)$ is the evaluation at $i$. 
%In other words, $\L X$ is a space of all continuous maps $l$ from the interval $[0,1]$
%to $X$ which satisfy $l(0)=\phi(l(1))$.

The twisted tube $\t_\phi X$ of
$X$ is defined by
\[
\t_\phi X = \frac{[0,1]\times X}{(0,x)\simeq (1,\phi(x))}
\]
and there is a canonical inclusion $\iota:X \hookrightarrow \t_\phi X$
defined by $\iota(x)=(0,x)$.
For example, $\t_{Id} X=S^1 \times X$ when $\phi$ is the identity map $Id$.
The cohomology of $\t_\phi X$ and $X$ are related by
the Wang exact sequence 
\begin{equation}\label{Wang}
  \cdots \to H^{n-1}(X) \xrightarrow{1-\phi^*} H^{n-1}(X) \xrightarrow{\delta}
  H^n(\t_\phi X) \xrightarrow{\iota^*} H^n(X) \xrightarrow{1-\phi^*} H^n(X) \to \cdots,
\end{equation}
where $\delta$ commutes with the action of $\a2$.
In particular, this exact sequence splits to short exact sequences
when $\phi^*$ is the identity map.

Let $ev: S^1 \times LX \to  X$ be the evaluation map defined by $ev(t,\gamma)=\gamma(t)$.
We define a map $\sigma_X:H^*(X) \to H^{*-1}(LX)$ by the following equation:
\[
  ev^*(x)= s \otimes \sigma_X (x) + 1 \otimes x, \quad (x \in H^*(X) ),
\]
where $s \in H^1(S^1)$ is the generator.
The \emph{twisted cohomology suspension} is defined by the following composition
\[
\hat{\sigma}_\phi :H^*(\t_\phi X) \xrightarrow{\sigma_{\t_\phi X}} H^{*-1}(L \t_\phi X))
 \xrightarrow{in^*} H^{*-1}(\L_\phi X),
\]
where $in:\L_\phi X \to L \t_\phi X$ is defined by $in(\gamma)(t)=(t,\gamma(t))$.
%\begin{eqnarray*}
% \L X &\to & L \t X \\
%   \gamma &\mapsto & t\mapsto (t,\gamma(t)).
%\end{eqnarray*}
When there exists a ring homomorphism $r: H^*(X)\to H^*(\t_\phi X)$ such that $\iota^*\circ r$ is the identity,
we define the corresponding \emph{module derivation} by $\sig = \hat{\sigma}_\phi \circ r: H^*(X) \to H^{*-1}(\L_\phi X)$.
We call such a homomorphism $r$ a section of $\iota^*$.
\begin{rem}\label{rem:LX}
  When $\phi$ is the identity map, we can take $r$ to be $\pi^*$,
   where $\pi:\t_{Id} X=S^1 \times X \to X$ is the projection.
   In this case, $\sig=\sigma_X$ coincides with Kuribayashi's module derivation
   $\mathcal{D}_X$~\cite{Kuribayashi}.
\end{rem}

The map $\sig$ together with the Wang
sequence above relates the cohomology of $X$ to that of $\L_\phi X$.
Consider the following conditions:
\begin{equation*}
  \label{eq:1} \left\{
  \begin{array}{cl}
    \mbox{(i)} & H^*(X) \mbox{ is a polynomial algebra } \z/2[x_1,x_2,\ldots,x_l],  \\
    \mbox{(ii)} &  \phi^* \mbox{ is the identity map.}
  \end{array} \right.
\end{equation*}
Then, the result of \cite{Kishimoto-Kono} specializes
to the following proposition.
\begin{prp}[Kishimoto-Kono]\label{KK-prop}
  Assume that the conditions (i) and (ii) are satisfied
  and there is a section $r$ of
  $\iota^*$ which commutes with the action of $\a2$.
  Denote by $e: \L_\phi X\to X$ the map defined by $e(\gamma)=\gamma(0)$.
  Then, we have
  \begin{enumerate}
  \item $\sig$ commutes with the action of $\a2$,
  \item $\sig(xy)= \sig(x)e^*(y) +
    e^*(x)\sig(y)$ for $x,y \in H^*(X)$,
  \item the elements 
  $\{\sig(x_1),\sig(x_2), \ldots ,\sig(x_l)\}$ form a simple system of generators for $H^*(\L_\phi X)$ 
    as an algebra over $\z/2[e^*(x_1),...,e^*(x_l)]$:
  \[
  H^*(\L_\phi X)\cong \z/2[e^*(x_1),e^*(x_2),\ldots,e^*(x_l)]\otimes
    \Delta (\sig(x_1),\sig(x_2), \ldots ,\sig(x_l)).\]
   \end{enumerate}
    This together with the action of $\a2$ on $H^*(X)$ determines the ring structure of $H^*(\L_\phi X)$ by 
    $(\sig(x_i))^2=\sig(Sq^d(x_i))$ with $d=|x_i|-1$.
   In particular, 
   we have an isomorphism $H^*(\L_\phi X)\cong H^*(L X)$ as algebras over $\a2$.
\end{prp}

 We consider the case when $X$ is the $2$-completion of either 
 $BSpin(7)$, $BSpin(8)$, $BSpin(9)$, $BF_4$, or $BDI(4)$.
 When $\phi$ is the identity map or an unstable Adams operation $\psi^q$,
 the conditions (i) and (ii) can be verified by a case-by-case analysis.
 We construct a section $r: H^*(BG)\to H^*(\t_{\psi^q} BG)$ which commutes with the Steenrod operations,
 and use the above proposition to compute $H^*(\L_{\psi^q} BG)\cong H^*(LBG)$.

%%%%%
\section{Computations for $G=Spin(7), Spin(8), Spin(9)$}
We recall the concrete presentations of
the mod $2$ cohomology of $BSpin(7),
BSpin(8)$ and $BSpin(9)$ with the action of $\a2$ from \cite{Quillen, Kono}.

$H^*(BSpin(7))\cong \z/2[w_4,w_6,w_7,w_8]$ and the action of $\a2$
is determined by
\[
\begin{array}{rcccc}
  & w_4 & w_6 & w_7 & w_8 \\
  Sq^1 & 0 & w_7 & 0 & 0 \\
  Sq^2 & w_6 & 0 & 0 & 0 \\
  Sq^4 & w_4^2 & w_4 w_6 & w_4 w_7 & w_4 w_8. \\
%  Sq^8 & 0 & 0 & 0 & w_8^2
\end{array}
\]
$H^*(BSpin(8))\cong \z/2[w_4,w_6,w_7,w_8,e_8]$ and the action of $\a2$ is determined by
\[
\begin{array}{rccccc}
  & w_4 & w_6 & w_7 & w_8 & e_8 \\
  Sq^1 & 0 & w_7 & 0 & 0 & 0 \\
  Sq^2 & w_6 & 0 & 0 & 0 & 0 \\
  Sq^4 & w_4^2 & w_4 w_6 & w_4 w_7 & w_4 w_8 & w_4 e_8. \\
%  Sq^8 & 0 & 0 & 0 & w_8^2 & e_8^2
\end{array}
\]
$H^*(BSpin(9))\cong \z/2[w_4,w_6,w_7,w_8,e_{16}]$ and the action of $\a2$
is determined by
\[
\begin{array}{rccccc}
  & w_4 & w_6 & w_7 & w_8 & e_{16} \\
  Sq^1 & 0 & w_7 & 0 & 0 & 0 \\
  Sq^2 & w_6 & 0 & 0 & 0 & 0 \\
  Sq^4 & w_4^2 & w_4 w_6 & w_4 w_7 & w_4 w_8 & 0 \\
  Sq^8 & 0 & 0 & 0 & w_8^2 & w_8 e_{16}+w_4^2 e_{16}. \\
%  Sq^{16} & 0 & 0 & 0 & 0 & e_{16}^2
\end{array}
\]

Based on these results,  we compute the mod $2$ cohomology of $LBG$ for $G=Spin(7), Spin(8)$, and $Spin(9)$.
Note that the conditions (i) is seen to be satisfied from the above concrete presentation.
The condition (ii) for $\phi=\psi^q$ is verified as follows.
%Recall that the rational cohomology $H^*(BG;\q)$ is isomorphic to the invariant ring $H^*(BT;\q)^W$, where $T$ is the maximal torus of $G$ and $W$ is its Weyl group.
The generators $w_4, w_8, e_8$, and $e_{16}$ are seen to be the mod $2$ reductions of torsion free integral classes by the Bockstein spectral sequence since the dimension of the $\z/2$-vector space $H^i(H^*(BSpin(n)), Sq^1)$ is same as the dimension of the $\q$-vector space $H^i(BSpin(n);\q)$ when $n=7,8,9$ and $i=4,8,16$. (In fact, Kono~\cite{bspin} showed that this is true for any $i$ and $n$ and the torsion elements of $H^*(BSpin(n);\z)$ are of order $2$.)
Since $(\psi^q)^*$ induces multiplication by a power of $q$ on the rational cohomology, it acts as the identity on these generators in the mod 2 cohomology.
By the compatibility of $(\psi^q)^*$ with the action of $\a2$, 
the other generators $w_6=Sq^2(w_4)$ and $w_7=Sq^1Sq^2w_4$ are also mapped identically by $(\psi^q)^*$
on the mod $2$ cohomology.

\begin{prp}\label{prp:bspin7}
  $H^*(LBSpin(7))\cong\z/2[v_4,v_6,v_7,v_8,y_3,y_5,y_7]/I \quad
  (|v_i|=i,|y_i|=i)$, where $I$ is the ideal generated by
   \[\{ 
    y_5^2 +y_3^2v_4+y_3v_7 , 
    y_3^4+ y_3^2 v_6 + y_5 v_7,
    y_7^2 +y_3^2v_8 +y_7 v_7\}.\]
    The action of  $\a2$
   is determined by
  \[
  \begin{array}{rccccccc}
    & v_4 & v_6 & v_7 & v_8 & y_3 & y_5 & y_7 \\
    Sq^1 & 0 & v_7 & 0 & 0       & 0 & y_3^2 & 0 \\
    Sq^2 & v_6 & 0 & 0 & 0       & y_5 & 0 & 0 \\ 
    Sq^4 & v_4^2 & v_4 v_6 & v_4 v_7 & v_4 v_8  & 0 & y_3v_6+y_5 v_4&
    y_3v_8+y_7 v_4. \\
%    Sq^8 & 0 & 0 & 0 & v_8^2 & 0 & 0 & 0
  \end{array}
  \]
  Moreover, $H^*(\L_{\psi^q} (BSpin(7))^{\wedge}_2)$ is isomorphic to $H^*(LBSpin(7))$ as algebras
  over $\a2$.
\end{prp}
\begin{proof}
When $\phi$ is the identity map, $\pi^*$ in Remark \ref{rem:LX}
serves as a section of $\iota^*$ which commutes with the action of $\a2$.
We apply Proposition \ref{KK-prop} 
  to obtain
   $H^*(LBSpin(7))\cong \z/2[v_4,v_6,v_7,v_8]\otimes
  \Delta[y_3,y_5,y_6,y_7]$,
  where $v_i=e^*(w_i)$ and $y_{i-1}=\sig(w_i) \quad (i=4,6,7,8)$.  
  The action of $\a2$ on $y_i \ (i=3,5,7)$ are determined as follows:
  \begin{eqnarray*}
    && Sq^1 y_3 = Sq^1 \sig(w_4) = \sig(Sq^1 w_4) =0 \\
    && Sq^2 y_3 = \sig(Sq^2 w_4) = \sig(w_6) =y_5 \\
    && Sq^1 y_5 = \sig(Sq^1 w_6) = \sig(w_7) =y_6 \\
    && Sq^2 y_5 = \sig(Sq^2 w_6) = 0 \\
    && Sq^4 y_5 = \sig(Sq^4 w_6) =
    \sig(w_4 w_6) = \sig(w_4) e^*(w_6) +\sig(w_6) e^*(w_4) 
    =   y_3 v_6 +y_5v_4 \\
    && Sq^1 y_7 = \sig(Sq^1 w_8) = 0 \\
    && Sq^2 y_7 = \sig(Sq^2 w_8) =0 \\
    && Sq^4 y_7 = \sig(Sq^4 w_8) = \sig(w_4 w_8)=y_3 v_8+ y_7 v_4. 
  \end{eqnarray*}
  With the aid of the Adem relations, we determine the ring structure
  as follows:
  \begin{align*}
     y_3^2 &= Sq^3 y_3 = Sq^1 Sq^2 y_3 = Sq^1 y_5 = y_6 \\
     y_5^2 &= Sq^5 y_5 = Sq^1 Sq^4 y_5 = Sq^1( y_3 v_6 + y_5 v_4)
    = y_3 v_7 +y_6 v_4 = y_3v_7 + y_3^2v_4\\
     y_7^2 &= Sq^7 y_7 = Sq^1 Sq^2 Sq^4 y_7 = Sq^1 Sq^2 (y_3 v_8+ y_7 v_4)=
     Sq^1(y_5 v_8 +y_7  v_6 )= y_3^2 v_8 +y_7 v_7 \\
      y_3^4 &=y_6^2=Sq^6 y_6=(Sq^2 Sq^4+Sq^5 Sq^1) y_6= \sig((Sq^2 Sq^4+Sq^5 Sq^1) w_7) \\
        &  = \sig(Sq^2 w_4 w_7)= \sig(w_6 w_7) = y_5 v_7 + y_3^2 v_6.
  \end{align*}

  To show that $H^*(\L_{\psi^q} (BSpin(7))^{\wedge}_2)$ is isomorphic to $H^*(LBSpin(7))$ as algebras
  over $\a2$ by Proposition \ref{KK-prop}, we have to construct a section
  $r: H^*(BSpin(7)^{\wedge}_2)\to H^*(\t_{\psi^q} (BSpin(7)))^{\wedge}_2)$ of $\iota^*$ in \eqref{Wang}
  which commutes with the action of $\a2$.
  To do so, we carefully choose an element
  $u_i \in (\iota^*)^{-1}(w_i) \subset H^i(\t_{\psi^q} (BSpin(7))^{\wedge}_2)$
  for each generator $w_i$ of $H^i(BSpin(7)^{\wedge}_2)$
   so that the action of $\a2$ on $u_i$ is compatible with that on $w_i$.

  The Wang sequences \eqref{Wang}
  give the following commutative diagram of exact sequences
  \[\xymatrix{
   {} \ar[r]^(0.15){1-(\psi^q)^*} & H^{n-1}(BSpin(7)^{\wedge}_2;\z/4) \ar[r]^{\delta} \ar[d]^\rho & H^{n}(\t_{\psi^q} (BSpin(7))^{\wedge}_2;\z/4)
  \ar[r]^{\iota^*} \ar[d]^\rho  & H^{n}(BSpin(7)^{\wedge}_2;\z/4) \ar[r]^(0.85){1-(\psi^q)^*} \ar[d]^\rho  & {} \\
  0 \ar[r] & H^{n-1}(BSpin(7)^{\wedge}_2;\z/2) \ar[r]^\delta & H^n(\t_{\psi^q} (BSpin(7))^{\wedge}_2;\z/2)
  \ar[r]^{\iota^*} & H^n(BSpin(7)^{\wedge}_2;\z/2) \ar[r] & 0,
  } \]
   where $\rho$ is the map in the
   Bockstein exact sequence
   \begin{align*}
   \to & H^n(\t_{\psi^q} (BSpin(7))^{\wedge}_2;\z/2) \to H^n(\t_{\psi^q} (BSpin(7))^{\wedge}_2;\z/4) \xrightarrow{\rho} 
   H^n(\t_{\psi^q} (BSpin(7))^{\wedge}_2;\z/2)  \\ 
   \xrightarrow{Sq^1} & H^{n+1}(\t_{\psi^q} (BSpin(7))^{\wedge}_2;\z/2) \to \cdots.
	\end{align*}   
  Since $w_4$ and $w_8$ are the mod $2$ reduction of torsion free integral classes,
  there are $\tilde{w}_i \ (i=4,8)$ in $H^*(BSpin(7)^{\wedge}_2;\z/4)$ such that $\rho(\tilde{w}_i)=w_i$. 
  Furthermore, $H^*(\psi^q;\z/4)$ acts as multiplication by one on $\tilde{w}_i \ (i=4,8)$.
  Hence, we can take $\hat{w}_4, \hat{w}_8\in H^{n}(\t_{\psi^q} (BSpin(7)))^{\wedge}_2;\z/4)$
  such that $\rho \iota^*(\hat{w}_i)=\iota^*\rho(\hat{w}_i)=w_i \ (i=4,8)$.
  Define $u_4=\rho(\hat{w}_4)$, $u_8=\rho(\hat{w}_8)$, $u_6=Sq^2 u_4$, and $u_7=Sq^1 u_6$.
   Since $Sq^1\circ \rho=0$, we have $Sq^1(u_i)=0$ for $i=4,8$.
   We compute $Sq^2(u_6)=Sq^2Sq^2(u_4)=0$.
   Since $\iota^*(Sq^4(u_6))=Sq^4\iota^*(u_6)=Sq^4(w_6)=w_4w_6=\iota^*(u_4u_6)$,
   $Sq^4(u_6)+u_4u_6$ is in the image of $\delta$. As $H^9(BSpin(7)^{\wedge}_2)=0$, we have $Sq^4(u_6)=u_4u_6$.
   Similarly, since $\iota^*(Sq^2 u_8)= Sq^2 \iota^*(u_8)=  Sq^2 w_8=0$, 
   we have $Sq^2 u_8=0$. 
    
   Since $H^{11}(BSpin(7)^{\wedge}_2) \cong \z/2$ is generated
   by $w_4 w_7$ and
   $\iota^*(Sq^4 u_8)= Sq^4 w_8= w_4 w_8=\iota^*(u_4u_8)$,
   we have $Sq^4 u_8=u_4 u_8+\epsilon \delta(w_4 w_7)$, where $\epsilon=0$ or $1$.
   To see $\epsilon=0$, we compute
   \[
   Sq^4 Sq^4 u_8 = Sq^4(u_4 u_8+\epsilon \delta(w_4 w_7)) = Sq^4(u_4 u_8) = u_4\epsilon \delta(w_4 w_7).
   \]
   Since $Sq^4 Sq^4 u_8=(Sq^6 Sq^2 + Sq^7 Sq^1)u_8=0$ by the Adem relation, 
   we have $\epsilon=0$.
 
  Take $r$ to be the ring homomorphism defined
  by $r(w_i)=u_i \ (i=4,6,7,8)$, then $r$
  is a section of $\iota^*$ which commutes with the action of $\a2$.
\end{proof}

The computation for $G=Spin(8)$ is similarly conducted.
\begin{prp}\label{prp:bspin8}
  $H^*(LBSpin(8))=\z/2[v_4,v_6,v_7,v_8,f_8,y_3,y_5,y_7,z_7]/I$
  $(|v_i|=i,|y_i|=i,|f_8|=8,|z_7|=7)$, where $I$ is the ideal generated
  by 
  \[ \{ y_5^2 +y_3^2 v_4+y_3v_7 , 
  y_3^4+ y_3^2v_6 + y_5 v_7, 
  y_7^2 +y_3^2 v_8 +y_7 v_7, 
  z_7^2 + y_3^2 f_8 +z_7 v_7 \}. \]
    The action of $\a2$ is determined by
  \[
  \begin{array}{rccccccccc}
    & v_4 & v_6 & v_7 & v_8 & f_8 & y_3 & y_5 & y_7 & z_7 \\
    Sq^1 & 0 & v_7 & 0 & 0 & 0      & 0 & y_3^2 & 0 & 0 \\
    Sq^2 & v_6 & 0 & 0 & 0 & 0       & y_5 & 0 & 0 & 0\\ 
    Sq^4 & v_4^2 & v_4 v_6 & v_4 v_7 & v_4 v_8 & v_4 f_8    & 0 & y_3v_6+
    y_5 v_4&
    y_3v_8+y_7 v_4 & y_3 f_8+z_7 v_4 . \\
%    Sq^8 & 0 & 0 & 0 & v_8^2 & f_8^2          & 0 & 0 & 0 & 0
  \end{array}
  \]
  Moreover, 
  $H^*(\L_{\psi^q} (BSpin(8))^{\wedge}_2)$ is isomorphic to $H^*(LBSpin(8))$ as algebras
  over $\a2$.
\end{prp}

For the case when $G=Spin(9)$, 
we make use of the following lemma.
\begin{lem}[Lahtinen]\label{lem:tube}
Assume $H^*(X)\cong \z/2[x_1,\ldots,x_l]$ and $\phi^*$ is the identity.
For $1\le i \le l$, let $u_i\in H^*(\t_\phi X)$ be any element 
such that $\iota^*(u_i)=x_i$.
Denote by $s\in H^1(\t_\phi X)$ the pullback of the generator $H^1(S^1)$ under the projection $\t_\phi X\to S^1$.
Then, we have
\[
H^*(\t_\phi X)\cong \z/2[u_1,\ldots,u_l,s]/(s^2)
\]
and $\delta(p(x_1,\ldots,x_l))=p(u_1,\ldots,u_l)s$ for any polynomial $p\in \z/2[x_1,\ldots,x_l]$.
\end{lem}

\begin{prp}\label{prop:bspin9}
  $H^*(LBSpin(9))=\z/2[v_4,v_6,v_7,v_8,f_{16},y_3,y_5,y_7,z_{15}]/I$
  $ (|v_i|=i,|y_i|=i,|f_{16}|=16,|z_{16}|=16)$, where $I$ is the ideal
  generated by
  \[ \{
  y_5^2 +y_3v_7 +v_4y_3^2, 
  y_3^4+ y_3^2 v_6 + y_5 v_7,
  y_7^2 +y_3^2 v_8 +y_7 v_7, 
  z_{15}^2 + v_7 v_8 z_{15}  + v_7 y_7f_{16} +  y_3^2 v_8 f_{16}
  \}. \] 
  The action of $\a2$ is determined by
  \[
  \begin{array}{rccccccccc}
    & v_4 & v_6 & v_7 & v_8 & f_{16} & y_3 & y_5 & y_7 & z_{15} \\
    Sq^1 & 0 & v_7 & 0 & 0 & 0      & 0 & y_3^2 & 0 & 0\\
    Sq^2 & v_6 & 0 & 0 & 0 & 0       & y_5 & 0 & 0 & 0\\ 
    Sq^4 & v_4^2 & v_4 v_6 & v_4 v_7 & v_4 v_8 & 0    & 0 & y_3v_6+y_5v_4 &
    y_3v_8+y_7 v_4 & 0\\
    Sq^8 & 0 & 0 & 0 & v_8^2 & v_8 f_{16}  +v_4^2   f_{16}     
       & 0 & 0 & 0 &
    y_7 f_{16} +  v_8  z_{15}+   v_4^2 z_{15}.
  \end{array}
  \]
  Moreover, 
  $H^*(\L_{\psi^q} (BSpin(9))^{\wedge}_2)$ is isomorphic to $H^*(LBSpin(9))$ as algebras
  over $\a2$.
\end{prp}
\begin{proof}
  In dimensions up to $8$, the calculation of $H^*(LBSpin(9))$ is completely same as in
  the case of $H^*(LBSpin(7))$. 
  Let $\phi$ be the identity map and define $f_{16}=e^*(e_{16})$ and $z_{15}=\sig(e_{16})$.
  We have
  \begin{align*}
   Sq^8 z_{15} &=\sig(Sq^8 e_{16})=
   \sig(w_8 e_{16}+w_4^2 e_{16})=
y_7 f_{16}+v_8 z_{15}+ v_4^2 z_{15} \\
  z_{15}^2 &=
  Sq^{15}z_{15}=\sig (Sq^{15}e_{16}) = \sig (Sq^7Sq^8e_{16})
  =\sig(Sq^7 (w_8 e_{16}+w_4^2 e_{16})) \\
  &= \sig(Sq^3Sq^4(w_8)e_{16}+Sq^3Sq^4(w_4^2) e_{16})
  =\sig(Sq^1 Sq^2(w_4 w_8)e_{16}+ Sq^3(w_6^2)e_{16})\\
  &=  \sig(w_7 w_8 e_{16})= v_7v_8\sig(e_{16})+\sig(w_7w_8)f_{16}=
  v_7 v_8 z_{15}  + v_7 y_7f_{16} +  y_3^2 v_8 f_{16}.
  \end{align*}
  
  Now, we will construct a section $r: H^*(BSpin(9)^{\wedge}_2)\to H^*(\t_{\psi^q} (BSpin(9))^{\wedge}_2)$ of $\iota^*$ which commutes with the action of $\a2$.
  Define $u_i \ (i=4,6,7,8)$ as in the case of $BSpin(7)$.
  We can choose an element $h_{16}\in H^{16}(\t_{\psi^q} (BSpin(9))^{\wedge}_2)$
  such that $\iota^*(h_{16})=e_{16}$ and $Sq^1(h_{16})=0$
   by the same argument for $u_8$ as in the case of $BSpin(7)$.
   By the Wang sequence for $\z/2$, we have $Sq^2 h_{16} =
  \epsilon_1 \delta(w_4 w_6 w_7)$ since $H^{17}(BSpin(9)) \simeq
  \z/2$ is generated by $w_4 w_6 w_7$  and
  $\iota^*(Sq^2 h_{16})=Sq^2 e_{16}=0$.
  We see $Sq^2 Sq^2 h_{16}=\epsilon_1 \delta(Sq^2(w_4 w_6 w_7)) =\epsilon_1
  \delta(w_6^2 w_7)$.
  Since $Sq^2 Sq^2 h_{16} =Sq^3 Sq^1 h_{16} =0$ by the Adem relation, we have $\epsilon_1=0$.

  Similarly, we have 
  \[
  Sq^4 h_{16}= \epsilon_2 \delta(w_4^3 w_7) +
  \epsilon_3 \delta(w_6^2 w_7) + \epsilon_4 \delta(w_4 w_7 w_8)
  \]
  and
  \[
  Sq^4 Sq^4 h_{16}=(\epsilon_2+\epsilon_3) \delta(w_4 w_6^2 w_7) +
   \epsilon_4 \delta(w_4^2 w_7 w_8).
  \] 
  Since $Sq^4 Sq^4 h_{16}= (Sq^7 Sq^1 + Sq^6
  Sq^2)h_{16} =0$, we have $\epsilon_2=\epsilon_3, \epsilon_4=0$.  
  Put $h'_{16}=h_{16}+ \epsilon_2 \delta(w_4^2 w_7)$, then we
  have $Sq^4 h'_{16}=0$ since $Sq^4 (w_4^2 w_7)=w_4^3 w_7 + w_6^2 w_7$.
   Since $Sq^i (w_4^2 w_7)=0 \ (i=1,2)$, we have $Sq^i h'_{16}=0 \ (i=1,2)$.

  Again by the Wang sequence, we have 
  \[
  Sq^8 h'_{16}= u_8 h'_{16}+u_4^2 h'_{16}+ \epsilon_5
  \delta(w_4^4 w_7)+ \epsilon_6 \delta(w_4^2 w_7 w_8) + \epsilon_7
  \delta(w_4 w_6^2 w_7) + \epsilon_8 \delta(w_7 w_8^2) + \epsilon_9
  \delta(w_7 e_{16})
  \]
  and by Lemma \ref{lem:tube} we have
  \[
  Sq^8 Sq^8 h'_{16}=su_7\left(\epsilon_5 u_4^6 + (\epsilon_5+\epsilon_6)u_4^4u_8+
  \epsilon_7u_4^3u_6^2+\epsilon_7u_4u_6^2u_8+(\epsilon_5+\epsilon_7)u_6^4+\epsilon_8u_8^3\right).
  \]
  Since
  $Sq^8 Sq^8 h'_{16}=(Sq^{12}Sq^4+Sq^{14}Sq^2+Sq^{15}Sq^1)h'_{16}=0$, we have
  $\epsilon_5=\epsilon_6=\epsilon_7=\epsilon_8=0$.
  Put $u'_{8}=u_8+\epsilon_9 su_7$. 
  We see
  \begin{align*}
  Sq^8 h'_{16} =u_8 h'_{16}+u_4^2 h'_{16}+\epsilon_9 s u_7 h'_{16}= u'_8 h'_{16}+u_4^2 h'_{16}.
  \end{align*}	
  Now, a section $r$ of $\iota^*$ is obtained by setting
  $r(w_i)=u_i \ (i=4,6,7)$, $r(w_8)=u'_8$, and $r(e_{16})=h'_{16}$.
  By the proof of Proposition \ref{prp:bspin7},
  we know that any lift $u'_8$ of $w_8$ behaves similarly with respect to the action of $\a2$,
  and hence, $r$ commutes with the action of $\a2$.
  \end{proof}

%%%%%%
\section{Computations for  $G=F_4$}
Denote by $i$ the classifying map of the canonical inclusion
$Spin(9)\hookrightarrow F_4$, where $F_4$ is the compact exceptional Lie group of type $F_4$.
Kono~\cite{Kono} showed that $i^*$ is injective and the mod $2$ cohomology of $BF_4$ over $\a2$
was determined as follows:
\[ H^*(BF_4)= \z/2[x_4,x_6,x_7,x_{16},x_{24}], \]
where $i^*(x_4)=w_4,
i^*(x_6)=w_6, i^*(x_7)=w_7, i^*(x_{16})=e_{16}+w_8^2, i^*(x_{24})=w_8
e_{16}$,   
and the action of $\a2$ is determined by
\[
\begin{array}{rcccccc}
  & x_4 & x_6 & x_7 & x_{16} & x_{24} \\
  Sq^1 & 0 & x_7 & 0 & 0 & 0 \\
  Sq^2 & x_6 & 0 & 0 & 0 & 0 \\
  Sq^4 & x_4^2 & x_4 x_6 & x_4 x_7 & 0 & x_4 x_{24} \\
  Sq^8 & 0 & 0 & 0 & x_{24}+x_4^2 x_{16} & x_4^2 x_{24}\\
  Sq^{16} & 0 & 0 & 0 & x_{16}^2 & x_{16}x_{24}+x_4 x_6^2 x_{24}.
\end{array}
\]
The condition (ii) for $\phi=\psi^q$ is satisfied as $i^*$ is compatible with $\psi^q$.

\begin{prp}\label{prp:bf4}
  $H^*(LBF_4)=\z/2[v_4,v_6,v_7,v_{16},v_{24},y_3,y_5,y_{15},y_{23}]/I$
  $(|v_i|=i,|y_i|=i)$, where $I$ is the ideal generated by
  \[ \{
   y_5^2 +y_3v_7 +v_4y_3^2,
   y_3^4+v_6 y_3^2 + y_5 v_7,
   y_{15}^2 + v_7 y_{23}+v_{24}y_3^2,
  y_{23}^2 + y_3^2v_{16}v_{24}+v_7v_{24}y_{15}+v_7 v_{16}y_{23} \}. \]
  The action of $\a2$ is determined by
  \[
  \begin{array}{rccccc}
    & v_4 & v_6 & v_7 & v_{16} & v_{24} \\
    Sq^1 & 0 & v_7 & 0 & 0 & 0 \\
    Sq^2 & v_6 & 0 & 0 & 0 & 0 \\
    Sq^4 & v_4^2 & v_4 v_6 & v_4 v_7 & 0 & v_4 v_{24} \\
    Sq^8 & 0 & 0 & 0 & v_{24}+v_4^2 v_{16} & v_4^2 v_{24} \\
    Sq^{16} & 0 & 0 & 0 & v_{16}^2 & v_{16}v_{24}+v_4 v_6^2 v_{24} \\
    \\
    \hline \\ 
    & y_3 & y_5 & y_{15} & y_{23} \\
    Sq^1      & 0 & y_3^2 & 0 & 0\\
    Sq^2       & y_5 & 0 & 0 & 0\\ 
    Sq^4     & 0 & y_3v_6+v_4y_5 &  0  &  y_3 v_{24} + v_4 y_{23}\\
    Sq^8   & 0 & 0 &  y_{23}+ v_4^2 y_{15} & v_4^2 y_{23}\\
    Sq^{16} & 0 & 0 & 0 & J,
  \end{array}
  \] where 
  $J=v_{24}y_{15}+v_{16}y_{23}+y_3 v_6^2 v_{24}+v_4 v_6^2 y_{23}$.
  Moreover, $H^*(\L_{\psi^q} (BF_4)^{\wedge}_2)$ is isomorphic to $H^*(LBF_4)$
   as algebras over $\a2$.
\end{prp}
\begin{proof}
  Let $\phi$ be the identity map, and put $v_i=e^*(x_i)$ and $y_{i-1}=\sig(x_i) \quad (i=4,6,7,16,24)$.  
  In dimensions up to $8$, the calculation of $H^*(LBF_4)$ is completely parallel to the
  the case of $H^*(LBSpin(9))$.
  The rest is computed as follows:
  \begin{align*}
    Sq^1 y_{15} &= \sig(Sq^1 x_{16}) =0 \\
     Sq^2 y_{15} &= \sig(Sq^2 x_{16}) =0 \\
     Sq^4 y_{15} &= \sig(Sq^4 x_{16}) =0 \\
     Sq^8 y_{15} &= \sig(Sq^8 x_{16}) =\sig(x_{24}+x_4^2x_{16})=y_{23}+v_4^2y_{15}  \\
     Sq^1 y_{23} &= \sig(Sq^1 x_{24}) =0 \\
     Sq^2 y_{23} &= \sig(Sq^2 x_{24}) =0 \\
     Sq^4 y_{23} &= \sig(Sq^4 x_{24}) =\sig(x_4 x_{24})=y_3 v_{24}+v_4 y_{23} \\
     Sq^8 y_{23} &= \sig(Sq^8 x_{24}) =\sig(x_4^2 x_{24})=v_4^2 y_{23} \\
     Sq^{16} y_{23} &= \sig(Sq^{16} x_{24}) =\sig(x_{16}x_{24}+x_4 x_6^2 x_{24})
    =y_{15}v_{24}+v_{16}y_{23}+y_3 v_6^2 v_{24}+v_4 v_6^2 y_{23} \\
     y_{15}^2&=Sq^{15}y_{15}=\sig(Sq^{15}x_{16})=\sig(Sq^7 Sq^8 x_{16})=\sig(Sq^7 x_{24})=
    \sig(x_7x_{24})=v_7 y_{23}+y_3^2 v_{24} \\
     y_{23}^2&=Sq^{23}y_{23}=\sig(Sq^{23}x_{24})=\sig(Sq^7 Sq^{16}x_{24}) 
     =\sig(Sq^7( x_{16}x_{24} +x_4 x_6^2 x_{24})) \\
    & \sig(x_{16}Sq^7(x_{24}) + Sq^7Sq^4(x_6^2 x_{24}))) = \sig(x_7x_{16}x_{24})=
    y_3^2v_{16}v_{24}+v_7y_{15}v_{24}+v_7 v_{16}y_{23}.
  \end{align*}
  
  Now, we will construct a section $r: H^*((BF_4)^{\wedge}_2)\to H^*(\t_{\psi^q} (BF_4)^{\wedge}_2)$ of $\iota^*$ which commutes with the action of $\a2$.
  By the uniqueness of the Adams operations~\cite{selfmap}, the following diagram is homotopy commutative:
  \[ \xymatrix{
    BSpin(9)^{\wedge}_2 \ar[r]^{\psi^q_{BSpin(9)}}
     \ar[d]^{(i)^{\wedge}_2} & BSpin(9)^{\wedge}_2 \ar[d]^{(i)^{\wedge}_2} \\
    (BF_4)^{\wedge}_2 \ar[r]^{\psi^q_{BF_4}} & (BF_4)^{\wedge}_2 }. \]
  By the
  naturality of the construction of the twisted tube, 
  there is a map $\t_{\psi^q} (i)^{\wedge}_2:\t_{\psi^q} BSpin(9)^{\wedge}_2 
  \to \t_{\psi^q} (BF_4)^{\wedge}_2 $, which makes the following
  diagram commute:
  \[ \xymatrix@=30pt{
    0 \ar[r] & H^{n-1}(BSpin(9)^{\wedge}_2) \ar[r]^{\delta_{BSpin(9)}} 
   & H^n(\t_{\psi^q} (BSpin(9))^{\wedge}_2) \ar[r]^{\iota^*_{BSpin(9)}} 
   & H^n(BSpin(9)^{\wedge}_2) \ar[r] & 0 \\
    0 \ar[r] & H^{n-1}((BF_4)^{\wedge}_2) \ar[r]^\delta \ar[u]^{i^*}
   & H^n(\t_{\psi^q} (BF_4)^{\wedge}_2) \ar[r]^{\iota^*} \ar[u]^{(\t_{\psi^q} (i)^{\wedge}_2)^*}  
   & H^n((BF_4)^{\wedge}_2) \ar[u]^{i^*} \ar[r] & 0,
  }
  \]
  where the horizontal lines are the Wang sequences and the vertical arrows are injections.
   Let $r_{BSpin(9)}$ be the section of $\iota^*_{BSpin(9)}$ constructed in the proof of 
   Proposition \ref{prop:bspin9}.
   Let $u_i=r_{BSpin(9)}(w_i) \ (i=4,6,7,8)$ and $h_{16}=r_{BSpin(9)}(e_{16})$.
 Since $i^*$ is an isomorphism in degrees up to $8$, so is 
 $(\t_{\psi^q} (i)^{\wedge}_2)^*$. We define 
\[
  r(x_i)= ((\t_{\psi^q}(i)^{\wedge}_2)^*)^{-1} u_i,
  \quad (i=4,6,7).
\]
If $r_{BSpin(9)}(i^*(x_{16}))=r_{BSpin(9)}(e_{16}+w_8^2)=h_{16}+u_8^2$ is in the image of $(\t_{\psi^q}(i)^{\wedge}_2)^*$,
we define $r(x_{16})=((\t_{\psi^q}(i)^{\wedge}_2)^*)^{-1}(h_{16}+u_8^2)$. 
Otherwise, we replace $h_{16}$ in the following manner.
Since $(\iota^*)^{-1}(e_{16}+w_8^2) = h_{16}+u_8^2 + \mathrm{Im}\delta_{BSpin(9)}$
and $H^{15}(BSpin(9))/H^{15}(BF_4) \cong \z/2\{w_7w_8\}$, 
the element $h_{16}+u_8^2+\delta(w_7w_8)$ should lie in the image of $(\t_{\psi^q}(i)^{\wedge}_2)^*$ in this case.
We define $r(x_{16})=((\t_{\psi^q}(i)^{\wedge}_2)^*)^{-1}(h_{16}+u_8^2+\delta(w_7w_8))$.
Finally, noting $x_{24}=Sq^8(x_{16})+x_4^2 x_{16}$,
we define $r(x_{24})=Sq^8(r(x_{16}))+r(x_4^2 x_{16})$.
From the proof of Proposition \ref{prop:bspin9}, we see that $r$ commutes with the action of $\a2$.

\end{proof}

\section{Computations for $G=DI(4)$}
In \cite{Dwyer-Wilkerson}, Dwyer and Wilkerson constructed a finite loop space $DI(4)$,
whose classifying space $BDI(4)$ has the mod $2$ cohomology isomorphic to the mod $2$ Dickson invariant of rank
$4$, that is, $H^*(BDI(4))\cong \z/2 [x_8, x_{12}, x_{14}, x_{15}]$, where
$|x_j|=j$.  The action of $\a2$ is determined by
\[
\begin{array}{rcccc}
  & x_8 & x_{12} & x_{14} & x_{15} \\
  Sq^1 & 0 & 0 & x_{15} & 0 \\
  Sq^2 & 0 & x_{14} & 0 & 0 \\
  Sq^4 & x_{12} & 0 & 0 & 0 \\
  Sq^8 & x_8^2 & x_8 x_{12} & x_8 x_{14} & x_8 x_{15}. \\
\end{array}
\]
Notbohm~\cite{BDI(4)} showed that 
there is an injection $H^*(BDI(4))\to H^*(BSpin(7))$ compatible with the unstable Adams operation $\psi^q$,
 and hence, the condition (ii) is satisfied for $\phi=\psi^q$.

Grbi\'c~\cite{Grbic} calculated 
$H^*(BSol(q))\cong H^*(\L_{\psi^q} BDI(4))$ over $\a2$ by the Eilenberg-Moore spectral sequence.
Kuribayashi~\cite{Kuribayashi}
calculated the mod 2 cohomology of
 $H^*(LBDI(4))$ over $\a2$.
Here, we reproduce their results by the same method as in the previous sections.
\begin{prp}\label{prp:bdi4}
  $H^*(LBDI(4))\cong \z/2[v_8,v_{12},v_{14},v_{15},y_7,y_{11},y_{13}]/I$
  $(|v_i|=i,|y_i|=i)$, where $I$ is the ideal generated by
  \[ \{ y_{11}^2 + y_7v_{15}+v_8y_7^2, y_{13}^2 + y_{11}v_{15}+v_{12}y_7^2,
  y_7^4 +y_{13}v_{15} +v_{14}y_7^2 \}. \]
   The
  action of $\a2$ is determined by
  \[
  \begin{array}{rccccccccc}
    & v_8 & v_{12} & v_{14} & v_{15} & y_7 & y_{11} & y_{13} \\
    Sq^1 & 0 & 0 & v_{15} & 0        & 0 & 0 &  y_7^2 \\
    Sq^2 & 0 & v_{14} & 0 & 0        & 0 & y_{13} & 0  \\
    Sq^4 & v_{12} & 0 & 0 & 0        & y_{11} & 0 & 0 & \\
    Sq^8 & v_8^2 & v_8 v_{12} & v_8 v_{14} & v_8 v_{15}
    & 0 & v_8y_{11} + y_7v_{12} &  v_8y_{13} + y_7v_{14}. \\
  \end{array}
  \]
  Moreover, $H^*(\L_{\psi^q} BDI(4))$ is isomorphic to $H^*(LBDI(4))$ as algebras
  over $\a2$.
\end{prp}
\begin{proof}
 Let $\phi$ be the identity map, and put $v_i=e^*(x_i), \ y_{i-1}=\sig(x_i) \quad (i=8,12,14,15)$.
Just as in the previous sections, we have
  \begin{eqnarray*}
    && y_7^2=Sq^7 y_7 = \sig(Sq^7 x_8)= \sig(Sq^1Sq^2Sq^4 x_8) = \sig(x_{15}) = y_{14} \\
    && Sq^1 y_i = \sig(Sq^1 x_{i+1}) =0 \ (i=7,11) \\
    && Sq^1 y_{13} = \sig(Sq^1 x_{14}) = \sig(x_{15})= y_{14} = y_7^2 \\
    && Sq^2 y_i = \sig(Sq^2 x_{i+1}) =0 \ (i=7,13) \\
    && Sq^2 y_{11} = \sig(Sq^2 x_{12})= \sig(x_{14})= y_{13} \\
    && Sq^4 y_i = \sig(Sq^4 x_{i+1}) =0 \ (i=11,13) \\
    && Sq^4 y_7 = \sig(Sq^4 x_8)= \sig(x_{12})= y_{11} \\
    && Sq^8 y_7 = \sig(Sq^8 x_8) =0  \\
    && Sq^8 y_{11} = \sig(Sq^8 x_{12})= \sig(x_8 x_{12})= y_7 v_{12}+ v_8y_{11} \\
    && Sq^8 y_{13} = \sig(Sq^8 x_{14})= \sig(x_8 x_{14})= y_7 v_{14}+ v_8y_{13} \\
    && y_{11}^2 = Sq^{11} y_{11}= \sig(Sq^{11} x_{12})=
    \sig(Sq^1 Sq^2 Sq^8 x_{12})= \sig(x_8 x_{15})= v_8 y_7^2+ y_7 v_{15} \\
    && y_{13}^2 = Sq^{13} y_{13} = \sig(Sq^{13}x_{14})=
    \sig((Sq^1 Sq^4 Sq^8 + Sq^{11} Sq^2)x_{14}) \\
    && \qquad =\sig(Sq^1Sq^4 x_8 x_{14})
    =\sig(x_{12} x_{15})= y_{11}v_{15}+v_{12}y_7^2 \\
    && y_7^4=y_{14}^2=Sq^{14}y_{14}=\sig(Sq^{14}x_{15})=
    \sig((Sq^{11}Sq^{3}+Sq^{2}Sq^{4}Sq^{8})x_{15})=\sig(x_{14}x_{15})
    = y_{13}v_{15}+v_{14}y_7^2.
  \end{eqnarray*}

We can choose an element $u_8\in H^{8}(\t_{\psi^q}BDI(4))$
  such that $\iota^*(u_8)=x_{8}$ and $Sq^1(u_{8})=Sq^1Sq^4(u_{8})=0$
   by the same argument for $u_8$ as in the proof of
   Proposition \ref{prp:bspin7}.  
    Put $u_{12}=Sq^4u_{12}, u_{14}=Sq^2 u_{12}$,
  and $u_{15}=Sq^1 u_{14}$.  Then, we have $Sq^1 u_i = 0 \
  (i=8,12,15)$. Since $H^{10}(\t_{\psi^q}BDI(4))=0$, we have $Sq^2 u_8 = 0$.
  Since $H^{19}(BDI(4))=0$ and $\iota^*(Sq^8 u_{12})=x_8 x_{12}$,
  we have $Sq^8 u_{12}= u_8 u_{12}$.
  Moreover, we have
  \begin{eqnarray*}
   && Sq^4 u_{12} = Sq^4 Sq^4 u_8= (Sq^7Sq^1+Sq^6Sq^2)u_8=  0 \\
    && Sq^2 u_{14}= Sq^2 Sq^2 u_{12}=0 \\
    && Sq^4 u_{14}= Sq^4 Sq^6 u_8= (Sq^8 Sq^2+Sq^{10}) u_8 =0 \\
    && Sq^8 u_{14}= Sq^8 Sq^2 u_{12} = (Sq^4 Sq^6 + Sq^2 Sq^8 + Sq^9 Sq^1)u_{12}
    = u_8 u_{14} \\
    && Sq^2 u_{15}= Sq^2 Sq^7 u_8 = (Sq^8 Sq^1 + Sq^9) u_8 =0 \\
    && Sq^4 u_{15}= Sq^4 Sq^7 u_8 = (Sq^9 Sq^2+Sq^{11}) u_8 = 0 \\
    && Sq^8 u_{15}= Sq^8 Sq^1 u_{14} = (Sq^9+Sq^2Sq^7)u_{14} 
    = Sq^1 Sq^8 u_{14}= u_8 u_{15}.
  \end{eqnarray*}
  Therefore,
  we can construct a section $r: H^*(BDI(4))\to H^*(\t_{\psi^q} BDI(4))$ of $\iota^*$ by $r(x_i)=u_i \ (i=8,12,14,15)$ so that $r$ commutes with $\a2$.
\end{proof}

% {\scshape
% \begin{flushright}
% \begin{tabular}{l}
% Nakanishi Printing Co., Ltd. \\
% Shimotachiuri-Ogawa-Higashi, \\
% Kamikyoku, Kyoto 602-8048, \\
% Japan \\
% {\upshape e-mail: infor@nacos.com}\\
% \\
% Nakanishi Printing Co., Ltd. \\
% Shimotachiuri-Ogawa-Higashi, \\
% Kamikyoku, Kyoto 602-8048, \\
% Japan \\
% {\upshape e-mail: infor@nacos.com}\\
% \end{tabular}
% \end{flushright}
% }

\end{document}